\documentclass[12pt]{iopart}
\usepackage{iopams}
\usepackage{latexsym} 
\usepackage{xcolor}

\newtheorem{prop}{Proposition}      
\newtheorem{thm}{Theorem}
\newtheorem{lemma}{Lemma}
\newtheorem{cor}{Corollary}
\newtheorem{defn}{Definition}

\newcommand{\epf}{\hfill   \mbox{$\Box $}}

\newcommand{\cf}{\mbox{${\mathbb C}$}}
\newcommand{\rf}{\mbox{${\mathbb R}$}}

\newcommand{\Mod}[1]{\ (\mathrm{mod}\ #1)}

\begin{document}

\title[Nonlinear Fourier transform and probability distributions]
{Nonlinear Fourier transform and probability distributions}

\author{Pavle Saksida}
\address{Faculty of Mathematics and Physics, 
University of Ljubljana, Jadranska 21, 
1000 Ljubljana, Slovenia \\
}

\ead{Pavle.Saksida@fmf.uni-lj.si}

\begin{abstract}
The paper describes some probabilistic and combinatorial aspects of nonlinear Fourier
transform associated with the AKNS-ZS problems. In the first of the two main results, we show that
a family of polytopes that appear in a power expansion of the nonlinear Fourier transform
is distributed according to the beta probability distribution. We establish this result by
studying an Euler type discretization of the nonlinear Fourier transform. This approach provides
our second main result, discovering a novel discrete probability distribution that approximates
the beta distribution. The numbers of alternating ordered partitions of an integer into
distinct parts are distributed according to our new distribution. Using another discretization, 
we also find a formula
for the values of alternating ordered partitions into non-distinct parts. We find a connection
between this  discretization and the multinomial
distribution.
\end{abstract}


\maketitle

\section{Introduction}

The nonlinear Fourier transform is the central object of the inverse scattering transform theory used to solve and analyze the integrable nonlinear partial differential
equations. 
A significant class of integrable equations are the equations
of the AKNS-ZS type. This class contains, e.g. sine-Gordon and nonlinear Schroedinger equations.
The pioneering work was done by Ablowitz, Kaup, Newel and Segur, \cite{AKNS1},
\cite{AKNS2}, and by
Zakharov and Shabat, \cite{ZS}.
In this paper, we shall consider the nonlinear Fourier transform ${\cal F}$
which appears in the study of the {\it periodic} AKNS-ZS problems.
To every well-behaved function $u(x) \colon [0, 1] \to
\cf$ it assigns the doubly infinite sequence $\{ {\cal F}[u](n)\}_{n \in \mathbb{Z}}$ of
$SU(2)$ matrices, given by ${\cal F}[u](n) = (- 1)^n\Phi (x = 1, n)$, where
$\Phi(x, n)$ is the solution of the linear initial value problem
\[
\Phi_x(x, n) = L(x, n) \cdot \Phi(x, n), \quad \Phi(0, n) = I.
\]
The coefficient matrix $L(x, n)$ is given by
\[
L(x, n) = \pmatrix{  \pi i \, n & u(x) \cr
- \overline{u(x)} & - \pi i \,n }.
\]
The transformation ${\cal F}$ can be thought of as a non-linearisation of the usual Fourier transformation. Namely, we have
\[
{\cal F}[u](n) = I + \pmatrix{ 0 & F[u](n) \cr
- F[\overline{u}](- n] & 0 } +
\sum_{d = 2}^{\infty} A_d[u](n),
\]
where $F$ is the linear Fourier transform and $u \to A_d[u]$ are matrix-valued nonlinear operators.
The amount of literature on various aspects of the inverse scattering method is
vast, so we shall only mention a few works in which the Fourier analysis aspect is more
pronounced. The foundational work was done by the originators mentioned above.
Nonlinear Fourier transforms of functions, defined on $\rf$ and $\rf^+$, were studied  by I. Gelfan'd,
A. Fokas and B. Pelloni in \cite{FoGe}, \cite{FoSu}, \cite{Pe}, and in their other works.
A different, but closely related transformation is described by T. Tao and C. Thiele in \cite{TTIntro}. Some aspects of the transformation, defined above, were studied in 
\cite{NLFT-1} and \cite{NLFT-2}.

Below, we shall consider ${\cal F}$, together with two of its discretizations.
Many authors studied discretizations of transformations  similar to ${\cal F}$, but usually acting
on the functions defined on $\rf$ or $\rf^+$, see e.g. \cite{YK1}, \cite{YK2}, \cite{Wahls2}.
M. Ablowitz and J. Ladik discovered a discretization that preserves the integrability of the AKNS-ZS systems, see
\cite{AL}.

In this paper we shall describe some  probabilistic and combinatorial  aspects of ${\cal F}$ which stem from its nonlinearity. Our first result concerns a certain set of polytopes. For every positive integer $d$ and every
$l \in [0, 1]$,  the polytope   $\widehat{D}_d(l)$ is given by
\[ \fl
\widehat{D}_d(l)  = \{(x_1, x_2, \ldots, x_d)  \in \rf^d; \  1 \geq x_1 \geq x_2 \geq \ldots \geq x_d \geq 0;
\ \ \sum _{j = 1}^d (- 1)^{j - 1} x_j = l\}.
\]
Denote by $D_d(l)$ the orthogonal projection of $\widehat{D}_d(l)$  on the hyperplane 
$\{(x_1, \ldots, x_{d - 1}, 0)\} \subset\rf ^d$.
These polytopes appear  in the power expansion of ${\cal F}[u_c](n)$ for the constant
function $u_c(x) \equiv u$.

We shall see that
\begin{equation}
{\cal F}[u_c](n) = I + \sum_{d = 1}^{\infty} u^d \int_0^1 {\rm Vol}(D_d(l)) 
\pmatrix{ 0 & e^{- 2 \pi i l n} \cr
               - e^{2 \pi i l n} & 0 }^d \,  dl.
\label{polytopes}
\end{equation}
\begin{thm}
For every dimension $d$, the volumes of polytopes $D_d(l)$ are essentially distributed according to the 
beta distribution with the shape parameters $(\frac{d}{2}, \frac{d}{2} + 1)$, if $d$ is even, and $(\frac{d +1}{2}, \frac{d + 1}{2})$, if $d$ is odd.
More concretely, we have the following expression:
\begin{equation} \fl
{\rm Vol}(D_d(l)) = \frac{1}{d!} \,  \left \{
\begin{array}{cl}  \frac{1}{B(\frac{d}{2},  \frac{d}{2} + 1)}  
\ l^{\frac{d}{2} - 1} (1 - l)^{\frac{d}{2}}   = p_{\beta}(l; \frac{d}{2}, \frac{d}{2} + 1) \ ; &
d \ \ {\rm even} \\
 &  \\
\frac{1}{B(\frac{d + 1}{2},  \frac{d + 1}{2})} \ 
  l^ {\frac{d - 1}{2} }  (1 - l)^{\frac{d - 1}{2}}
 = p_{\beta}(l; \frac{d + 1}{2}, \frac{d + 1}{2})
 \ ; & d \ \ {\rm odd}\ \ ,
 \end{array}
  \right.  
  \label{BetaVolume}
\end{equation}
where $p_{\beta}(l; a, b)$ denotes the probability density function of the distribution 
${\rm Beta}(a, b)$.
 \label{Beta1}
 \end{thm}
Expressions (\ref{polytopes}) and (\ref{BetaVolume})  point to 
the importance of the beta distribution
 for the nonlinear Fourier transform ${\cal F}$. In the proof of the above theorem, we shall use the relation between the
polytopes and ${\cal F}$ in an essential way.
The beta distribution is one of the oldest and most important probability distributions with a broad spectrum
of applications in different areas of probability and statistics, particularly in Bayesian
statistical inference. In recent times it is mentioned in virtually every book on machine learning and related topics. 
The beta distribution also appears connected with polytopes, although in a setting very different from ours. In
\cite{K-T-Z} and in many other works,
Kabluchko, Thale, and Zaporozhets, together with coworkers, describe exciting results
concerning the relations between volumes and
angles of random polytopes on the one hand, and beta distributions on the other.

To formulate and to  prove theorem \ref{Beta1} we used an appropriate discretization of ${\cal F}$,
namely the transformation ${\cal F}_N$ described in section \ref{Transforms}. The study of this discretization leads to another result
concerning beta distribution. 

\begin{thm}
Let the discrete probability distribution ${\rm Beta}_N(a, b)$ be given by the probability mass function
\[
P_N(\lambda; a, b) = c(N) \, \frac{(a + b + 1)!}{N^{a + b}} \ 
{ N \lambda - 1 \choose a }  { N - N \lambda \choose b}
\]
which is defined on the set of values $\lambda \in 
\{0, \frac{1}{N}, \frac{2}{N}, \ldots, \frac{N - 1}{N}\}$.
All the parameters, except for $\lambda$, are integers and
$C(N)$ is the normalizing factor.
This probability distribution is a discrete approximation of the continuous beta distribution
with the probability density function
\[
p_{\beta}(x; a, b) = \frac{1}{B(a + 1, b + 1)} \ x^{a} ( 1 - x)^{b}, 
\quad x \in [0, 1].
\]
More concretely, let the sequence of integers $\{l_N\}_{N \in \mathbb{N}}$ be such that  $\lambda_N < N$ and
$\lim _{N \to \infty} \frac{l_N}{N} = \lambda \in [0, 1]$. Then 
we have
\[
\lim _{N \to \infty} P_N(l_N; a, b) = p_{\beta}(\lambda; a, b).
\]
We also have
$
c(N) = 1 + {\cal O}(\frac{1}{N}).
$
\label{Beta2}
\end{thm}
A different approach to discretize the beta distribution is introduced by A. Punzo in \cite{APunzo}.

Theorem  \ref{Beta2} stems from the study of  the numbers $AQ_N(l, d)$
which count the {\it ordered alternating partitions of $l$ into $d$ distinct parts not greater than $N- 1$}, 
\begin{equation} \fl
AQ_N(l, d) = \sharp \{(l_1, l_2, \ldots, l_d); \ N - 1 \geq l_1 > l_2 > \ldots >  l_d \geq 0; \  \
\sum_{j = 1}^{d} (- 1)^{j - 1}  l_j = l \}.
\label{AQdef}
\end{equation}
We shall prove that the numbers $AQ_N(l, d)$ are essentially distributed according to the distribution $P_N$, 
given in theorem \ref{Beta2}. This fact is not surprising, since the numbers 
$AQ_N(l, d)$ can be viewed as discretizations of the volumes ${\rm Vol}(D_d(l))$. Just as 
the numbers ${\rm Vol}(D_d(l))$ are closely related to ${\cal F}$, so are the numbers
$AQ_N(j, d)$ closely related to ${\cal F}_N$. In proposition \ref{FN and AQ} of section 3 we prove the 
following formula which gives the expression of $AQ_N(l, d)$ in terms of the transformation 
${\cal F}_N$:
\[ \fl
AQ_N(l, d) = \frac{N^d}{d!} \ (\frac{d}{d u})^d|_{u = 0} \Bigl(\sum_{n = 0}^{N - 1}
\pmatrix{ e^{ 2 \pi i \frac{l n}{N}} & 0 \cr
                          0         &   e^{ - 2 \pi i \frac{l n}{N}} }\cdot {\cal F}_N[\vec{u}_c](n) 
\cdot  \pmatrix{0 &    1 \cr
               - 1 & 0 }^{- d} \,
\Bigr)_{1, 1},
\]
where $\vec{u}_c = (u, u, \ldots, u)$ is the constant discrete function.

Let now $AP_N(l, d)$ denote the number of ordered alternating partitions of $l$ into 
$d$ {\it non-distinct} parts not greater than $N- 1$. It turns out that finding the values
of $AP_N(l, d)$ demands a different approach. In section  \ref{Transforms}, we introduce another
discretization ${\cal G}_N$ of ${\cal F}$ and this transform enables us to find 
formulae for $AP_N(l, d)$.
In section \ref{partitions}, we prove the following proposition.
\begin{prop}
The number  $AP_N(l, d)$ of ordered alternating partitions of number $l$ into $d$ non-distinct
parts is given by 
\[ \fl
AP_N(l, d) =   \frac{N^d}{d!}\ {\cal D}_d \Bigl( \sum_{n = 0}^{N - 1} 
\pmatrix{ e^{ 2 \pi i \frac{l n}{N}} & 0 \cr
                          0         &   e^{ - 2 \pi i \frac{l n}{N}} }
D^d{\cal G}_N[\vec{u}](n)\cdot 
\pmatrix{0 &    1 \cr
               - 1 & 0 }^{- d} \,
\Bigr)_{1, 1} |_{\vec{u} = 0} . 
\]
Here $\vec{u} = (u_0, u_1, \ldots, u_{N - 1}) \in \rf^d$ and the operators ${\cal D}_d$ and $D^d$ are defined 
by
\[ \fl
{\cal D}_d = \sum_{\vec{k} \in [0, d]^N }
 \frac{\partial^{\sum_{j = 0}^{N - 1} k_j}}{\partial^{k_0} u_0 \,  \partial^{k_1}u_{k_1}\cdots \partial ^{k_{N - 1}}u_{N - 1}}
 \quad {\rm and} \quad
D^d(f(\vec{u})) = (\frac{d}{d s})^d|_{s = 0} \ f( s \vec{u}).
\]
\label{APexpression}
\end{prop}
The subscript $(1, 1)$ denotes the upper right term of the $2 \times 2$ matrix.

A central object in the study of the numbers $AQ_N(l, d)$ is the vector 
$ \vec{l} = (l_1, l_2, \ldots, l_d)$
together with its alternating sum $l = \sum_{j = 1}^d (- 1)^{j - 1} l_j$.
This vector appears in the study of $AP_N(l, d)$ in an implicit way. Let
$\vec{k} = (k_0, k_1,\ldots, k_{N - 1})$ be a vector of nonnegative integers. Then the analogue 
of $\vec{l}$ is
 the vector of those indices $(l_1, l_2, \ldots l_d)$ from $\{0, \ldots , N - 1\}$ for which
the components $k_{l_j}$ of $\vec{k}$ are {\it odd} integers. The analogue of the alternating sum
is  the function ${\rm alt}$, given by
$
{\rm alt} (\vec{k}) =  \sum_{j = d}^1 (- 1)^{d - j} l_j.
$
In section \ref{partitions}, we give a longer but clearer description of ${\rm alt}(\vec{k})$.

The discretization ${\cal G}_N$ turns out to be related to the multinomial distribution. The connection between the two objects can be seen in a variety of ways. One of them is the following proposition which we prove in section \ref{partitions}.
\begin{prop}
Let $\vec{X} = (X_0, X_2, \ldots, X_{N- 1})$ be a random vector with values in $(\mathbb{N} \cup \{0\})^N$,
and let the probability of the event $\vec{X} = \vec{k}$ be given by the multinomial distribution
\[
P(\vec{X} = \vec{k}) =  P_{N, d}(\vec{u}, \vec{k}) =  { d \choose k_0, k_1 \ldots k_{N - 1}} 
u_0^{k_0} u_1^{k_1} \cdots u_{N - 1}^{k_{N - 1}},
\]
where $\sum_{j = 0}^{N - 1} u_j = 1$.
Then the probability $P_{alt}(l)$ of the event that $\vec{X}$ will assume a value $\vec{k}$ with
$ {\rm alt} (\vec{k}) = l$ is equal to
\[ \fl
P_{alt}(l) = \sum_{{\rm alt}(\vec{k}) = l} P_{N, d}(\vec{u}, \vec{k})  =
\Bigl( N^d \Bigl( \sum_{n = 0}^{N - 1} 
\pmatrix{ e^{  2 \pi i \frac{l n}{N}} & 0 \cr
                             0       &     e^{- 2 \pi i \frac{l n}{N}} }
 \cdot D^d{\cal G}_N[\vec{u}](n) \Bigr) \cdot  \pmatrix{0 &    1 \cr
               - 1 & 0 }^{- d}  
\Bigr)_{1, 1},
\]
where ${\cal G}_N$ is the discretization of ${\cal F}$ appearing in proposition  \ref{APexpression}.
\label{multinomial}
\end{prop}

\section{Discrete nonlinear Fourier transforms}
\label{Transforms}

We have defined the nonlinear Fourier transform of functions $u(x) \colon [0, 1] \to \cf$
in the introduction. Definition in this form is usually given in the texts which study the integrable
ANKS-ZS equations. We shall rather represent ${\cal F}$ in a different
gauge. Let $G(x, n) = {\rm diag}( e^{ - \pi i n x}, e^{\pi i n x})$ be the (diagonal)  matrix of our
gauge transformation. The transformed coefficient matrix is
$
L^G(x, n) = G_x \cdot G^{ - 1}(x, n) + G(x, n) \cdot L(x, n) \cdot G^{ - 1}(x, n).
$
Its explicit expression is
\begin{equation}
L^G(x, n) = \pmatrix{ 0 & e^{ - 2 \pi i n x} u(x) \cr
- e^{ 2 \pi i n x} \overline{u(x)} & 0 }.
\label{LG}
\end{equation}
In the new gauge ${\cal F}[u](n)$ becomes ${\cal F}^G[u](n) = \Phi^G(x = 1, n)$, where
$\Phi^G(x, n)$ is the solution of the linear initial value problem
\begin{equation}
\Phi^G_x(x, n) = L^G(x, z) \cdot\Phi^G (x, n), \quad \Phi^G(0, n) = I.
\label{G-initial}
\end{equation}
We have $\Phi^G(x, n) = G(x, n) \cdot \Phi(x, n)$ and, since $n\in \mathbb{Z}$, we have
$
{\cal F}[u](n) = {\cal F}^G[u](n).
$

The solution to the problem (\ref{G-initial}) can be given in the form of the Dyson series.
\begin{equation}
\Phi^G(x, n) = I+ \sum_{ d = 1}^{\infty} \int_{\Delta_d(x)}
L^G(x_1, n) \cdot L^G(x_2, n) \cdots L^G(x_d, n) \ d\vec{x},
\label{Dyson}
\end{equation}
where $\Delta_d(x)$ is the ordered simplex of dimension $d$ with the edge length equal to $x$,
\[
\Delta_d(x) = \{(x_1, x_2, \ldots, x_d) \in \rf ^d; x \geq x_1 \geq x_2
\geq \ldots \geq x_d \geq 0 \}.
\]
Let us denote
\begin{equation}
E(x, n) = \pmatrix{ e^{\pi i x n} & 0 \cr
0 & e^{ - \pi i x n} }, \quad
J = \pmatrix{ 0 & 1 \cr
- 1 & 0 },
\label{notation E-J}
\end{equation}
and let $u(x)$ be real valued. Then we have
$
L^G(x, n) = u(x)\,E(- 2 x, n) \cdot J.
$
Matrices $E(x, n)$ and $J$ do not commute. Instead,  have the relation
\begin{equation}
E(x, n) \cdot J = J \cdot E(- x, n).
\label{commutationrelation}
\end{equation}
Using  (\ref{commutationrelation}) in the Dyson series and evaluating at $x = 1$ gives
\[ \fl
{\cal F}[u](n) = I + \sum_{d = 1}^{\infty}\int _{\Delta_d(1)} \! \! \! 
u(x_1) \, u(x_2) \cdots u(x_d) \, E\Bigl(- 2( \sum_{j = 1}^d (- 1)^{j - 1} x_j), n\Bigr)
\cdot J^d \
d\vec{x}
\]
which, upon setting $x_1 - x_2 + \ldots + (- 1)^{d - 1} x_d = l$, can be rewritten as
\begin{eqnarray} \fl
{\cal F}[u](n) & = & I + \sum_{d = 1}^{\infty} \int _0^1 E(- 2 l, n) \ \Bigl( \int_{\widehat{D}_d(l)} 
u(x_1) \, u(x_2) \cdots u(x_d) \, d_l \vec{x}  \Bigr)\ \cdot J^d   \frac{1}{\sqrt{d}} \,dl   \nonumber\\
 \fl & = & 
 I + \sum_{d = 1}^{\infty} \int _0^1 E(- 2 l, n) \ \Bigl( \int_{D_d(l)} 
\! \! \! \! \! {\cal U}(x_1, x_2, \ldots, x_{d - 1}; l)   \,d x_1 \cdots d x_{d - 1} \Bigr) 
J^d \, dl,
\label{Fexpansion1}
\end{eqnarray}
where the polytope $\widehat{D}_d(l)$ is given by
\[
\widehat{D}_d(l) = \{(x_1, x_2, \ldots x_d) \in \Delta_d(1);
\ \ \sum _{j = 1}^d (- 1)^{j - 1} x_j = l\},
\]
and $D_d(l)$ is its projection on the hyperplane $x_d = 0$. We denoted
\[ \fl
{\cal U}(x_1, x_2, \ldots, x_{d - 1}; l) )  = u(x_1) \cdots u(x_{d - 1})
u((- 1)^{d - 1}(l - (x_1 - x_2 + \ldots  + (- 1)^{d - 2} x_{d - 1})).
\]
By means of some linear algebra one can show that the volume 
forms $d_l \vec{x}$ and the Euclidean form $d x_1 \cdots d x_{d - 1}$ are related by
$d_l\vec{x} = \sqrt{d} \ d x_1 \cdots d x_{d - 1}$.

In the case where $u_c(x) \equiv u$
is a constant function, we get
\begin{equation}
{\cal F}[u_c](n) = I + \sum_{d = 1}^{\infty} u^d \int _0^1 {\rm Vol}(D_d(l)) \, E(- 2 l, n)
\cdot J^d  \,dl.
\label{expansion by volumes}
\end{equation}

\paragraph{Discretization ${\cal F}_N$}
We have obtained the nonlinear Fourier transform from
an initial value problem for a particular first-order linear differential equation.
An obvious approach to construct a discretization is to replace the differential equation with
a suitable difference equation. Let $\vec{u} = (u_0, u_1, \ldots, u_{N - 1}) \in \rf^N$ be
a vector which plays a role of a function of a discrete variable. Let the
$L$-matrix be given by
\[
L_N(k, n) = \pmatrix{ 0 & e^{ - 2 \pi i \frac{k n}{N}} u_k \cr
- e^{ 2 \pi i \frac{k n }{N}} u_k & 0 }.
\]
\begin{defn}
Let $k, n \in \{0, 1, \ldots, N - 1\}$.
Discrete nonlinear Fourier transform ${\cal F}_N[\vec{u}]$ of $\vec{u}$ is defined by
${\cal F}_N[\vec{u}](n) = \Phi_N(k = N - 1, n)$, where $\Phi_N$ is the solution
of the difference initial value problem
\[
N \Bigl( \Phi_N(k + 1, n) - \Phi_N(k, n) \Bigr) = L_N(k, n) \cdot \Phi_N(k, n), \quad
\Phi_N(0, n) =I.
\]
\end{defn}
Solving the above initial value problem and evaluating at $k = N - 1$ gives
\[
{\cal F}_N[\vec{u}](n) = \prod_{k = N - 1}^0 \Bigl(I + \frac{1}{N} L_N(k, n) \Bigr),
\]
and this can be expanded into
\begin{equation} \fl
{\cal F}_N[\vec{u}](n) = I + \sum_{ d = 1}^{N} \frac{1}{N^d}
\sum_{N - 1 \geq l_1> l_2 > \ldots > l_d \geq 0} \! \! \! \!
L_N(l_1, n) \cdot L_N(l_2, n) \cdots L_N(l_d, n).
\label{DiscDyson}
\end{equation}
This expression is a discrete analogue of the Dyson expansion (\ref{Dyson}).

Let us introduce the notation
\[
E_{\delta}(l, n) = E(l, \frac{n}{N}), \quad l, n \in \{0, 1, \ldots , N - 1\}
\]
where $E$ is given by (\ref{notation E-J}), and the subscript $\delta$ refers to the use in the discretized
context.
The coefficient matrix $L_N$ can be written in the form
\[
L_N(l, n) = u_l \ E_{\delta}(- 2 l, n) \cdot J,
\]
with $J$ also defined in (\ref{notation E-J}). By means of relation (\ref{commutationrelation}),
we can collect all the copies of $J$ in (\ref{DiscDyson}) on the right. Let $\vec{u}_c = (u, \ldots , u)$
be a constant vector. We get
\[ \fl
{\cal F}_N[\vec{u}_c](n) = I + \sum_{d = 1}^N (\frac{u}{N})^d
\sum_{N - 1 \geq l_1 > l_2 > \ldots > l_d \geq 0} \! \! \! \!
E_{\delta}\Bigl(- 2(l_1 - l_2 + \ldots + (- 1)^{d - 1} l_d), n\Bigr) \cdot J^d.
\]
If we denote $l = l_1 - l_2 + \ldots + (- 1)^{d - 1} l_d$, we can finally write
\begin{equation}
{\cal F}_N[\vec{u}_c](n) = I + \sum_{d = 1}^{N - 1}
(\frac{u}{N})^d \ \
\sum_{l = 0}^{N - 1} E_{\delta}(- 2 l , n) \sum_{(l_1, \ldots, l_d) \in 
\widehat{D}_d^{disc}(l)} J^d,
\label{FNexpansion}
\end{equation}
where
\begin{equation}
 \fl
\widehat{D}_d^{disc}(l) = 
\{(l_1, l_2, \ldots , l_d); \ N - 1 \geq l_1 > l_2 > \ldots > l_d \geq 0, \
\sum_{j = 1}^d (- 1)^{j - 1} l_j = l \}.
\label{Ddisc}
\end{equation}

\paragraph{Discretization ${\cal G}_N$} Another approach to discretize ${\cal F}$ is to start with
${\cal F}[u_s]$, where $u_s$ is a step function, $u_s(x) = \sum_{l = 0}^{N - 1} u_l \
\chi _{[l, l + (1/N)]}(x)$. This can be computed directly. We have
\[ \fl
{\cal F}[u_s](n) = {\rm Exp}(\frac{1}{N}L(N - 1, n)) \cdot {\rm Exp}(\frac{1}{N}L(N - 2, n))
\cdots {\rm Exp}(\frac{1}{N}L(0, n)).
\]
This discretization has its merits, but we will simplify it by separating the spatial and spectral
parameters:
\[
\frac{1}{N} L(l, n)= \pmatrix{\frac{\pi i}{N}\, n & 0 \cr
0 &  - \frac{\pi i}{N} \, n } +
\pmatrix{ 0 & \frac{u_l}{N} \cr
- \frac{u_l}{N} & 0 }.
\]
The Baker-Campbell-Hausdorff formula gives
\begin{eqnarray*} 
{\rm Exp}( \frac{1}{N}L(l, n))  & = &  \pmatrix{ e^{ \frac{\pi i}{N} n} & 0 \cr
0 & e^{ - \frac{\pi i}{N} n } } \cdot
\pmatrix{ \cos{\frac{u_l}{N}} & \sin{\frac{u_l}{N}} \cr
- \sin{\frac{u_l}{N}} & \cos{\frac{u_l}{N}} } + {\cal O}(\frac{1}{N^2})    \\ 
&  &  \\
&  = &
 E_{\delta}( 1 , n ) \cdot R(\frac{u_l}{N}) + {\cal O}(\frac{1}{N^2}).
\end{eqnarray*}
\begin{defn}
Let $l, n, \in \{0, 1, \ldots, N - 1\}$. Discrete nonlinear transform ${\cal G}_N[\vec{u}]$ of $\vec{u}$
is defined by
\[
{\cal G}_N[\vec{u}](n)= E_{\delta}(1, n) \cdot R(\frac{u_{N - 1}}{N}) \cdots
E_{\delta}( 1, n) \cdot R(\frac{u_1}{N}) \cdot E_{\delta}(1, n) \cdot R(\frac{u_0}{N}).
\]
\end{defn}
Using the obvious identity
$
E_{\delta}(1, n) = E_{\delta}(l + 1, n) \cdot E_{\delta}(- l, n)
$
we can express ${\cal G}_N$ in the form
\[ \fl
{\cal G}_N[\vec{u}](n) = \prod _{l = N - 1}^0 {\rm Ad}_{E_{\delta}(- l, n)}R(\frac{u_l}{N}) =
\prod _{l = N - 1}^0 \Bigl( \cos{\frac{u_l}{N}} \, I + 
\sin{\frac{u_l}{N}} \ E_{\delta}(- 2l , n) \cdot J \Bigr).
\]
It is easy to see that ${\cal F}_N$ and ${\cal G}_N$ are related by the formula
\[
{\cal G}_N[\vec{u}](n) = {\cal C}[\vec{u}] \ {\cal F}_N[\tan{(\vec{u}/N})](n),
\]
where ${\cal C}[\vec{u}] = \prod_{l = 0}^{N - 1} \cos{(u_l}/N)$, and
$\tan{(\vec{u}/N)} = (\tan{(u_0/N)}, \ldots, \tan{(u_{N - 1}/N)})$. We see that for small $\vec{u}$
the two discretizations differ only very little. We note that, unlike ${\cal F}_N$, the discretization
${\cal G}_N$ takes values in $SU(2)$.

\section{Ordered   alternating  partitions with distinct parts}

This section will first describe the connection between the numbers of ordered alternating
partitions with distinct parts $AQ_N(l, d)$ and the transformation ${\cal F}_N$. The results will lead us to the unexpected connection between the transformation ${\cal F}$, volumes of polytopes $D_d(l)$, and the
beta distribution, described in theorem \ref{Beta1}. We shall also provide a novel discretization
of the beta distribution, given by theorem \ref{Beta2}.

Let us recall definition (\ref{AQdef})
\[ \fl
AQ_N(l, d) = \sharp \{(l_1, l_2, \ldots, l_d); \ N - 1 \geq l_1 > l_2 > \ldots > l_d \geq 0; \ \
\sum_{j = 1}^{d} (- 1)^{j - 1} l_j = l \}.
\]
\begin{prop}
The power series expansion of ${\cal F}_N[u_c]$ around $u = 0$ is given by
\begin{equation}
{\cal F}_N[u_c](n) =
I + \sum_{d = 1}^N (\frac{u}{N})^d \ \sum_{ l = 0}^{N - 1} AQ_N(l, d) \, E_{\delta}( - 2 l, n) \cdot J^d .
\label{FN-AQ}
\end{equation}
The number $AQ_N(l, d)$ of alternating partitions of $l$ into $d$ distinct parts not greater than
$N - 1$ is given by the equation
\begin{equation}
AQ_N(l, d) = \frac{N^d}{d!} \ (\frac{d}{d u})^d|_{u = 0} \Bigl(\sum_{n = 0}^{N - 1}
E_{\delta}(2 l , n) \cdot {\cal F}_N[u_c](n) \cdot J^{ - d} \Bigr)_{1, 1}.
\label{AQ-FN}
\end{equation}
\label{FN and AQ}
\end{prop}

\noindent{\bf Proof:} The first formula of proposition follows immediately from equation (\ref{FNexpansion}). We only have to notice that $AQ_N(l, d)$ is equal to the
number of elements in $\widehat{D}_d^{disc}(l)$. To get (\ref{AQ-FN}), we multiply
both sides of (\ref{FN-AQ}) by $J^{- d}$ and then perform the inverse linear discrete
Fourier transform on both sides. To isolate the term containing the $d$-th power of $u$, we
take the $d$-th derivative with respect to $u$, evaluate at $u = 0$, and get   formula 
(\ref{AQ-FN}).

\epf

There is an explicit formula for the function $AQ_N(l, d)$.  We have:

\begin{prop}
For any $N \in \mathbb{N}$, $d \leq N$ and $l \in \{0, \ldots, N - 1\}$, we have
\begin{equation}
AQ_N(l, d) =
\left \{\begin{array}{cl}{ l - 1 \choose \lfloor \frac{d - 1}{2}\rfloor}
{N - l \choose \lfloor \frac{d}{2} \rfloor} \ ;
& d \ \ {\rm even} \\
& \\
{ l \choose \lfloor \frac{d - 1}{2} \rfloor} {N - l - 1\choose \lfloor \frac{d}{2}\rfloor} \ ; & d \ \
{\rm odd} \ \ .
\end{array}
\right.
\label{AQexplicit}
\end{equation}
\label{AQ-formula}
\end{prop}
Above we use the definition of the binomial symbol for which ${ a \choose b} = 0$ for negative $a$.

\noindent{\bf Proof:}
Let us define
\[  \fl
\widehat{AQ}_N(l, d) = \sharp \{ (l_1, \ldots , l_d); N \geq l_1 > \ldots > l_d \geq 1,
\ \ {\rm and} \ \ \sum_{j = 1}^d (- 1)^{j - 1} l_j = l \}.
\]
We claim that for $\widehat{AQ}_N(l, d)$ we have
\begin{equation}
\widehat{AQ}_N(l, d) = { l - 1 \choose \lfloor \frac{d - 1}{2}\rfloor}
{N - l \choose \lfloor \frac{d}{2} \rfloor}.
\label{altp-formula}
\end{equation}
The formula can be proved by induction on $N$. For $N = 2$, formula (\ref{altp-formula})
can be checked by hand. It is an easy exercise to show that $\widehat{AQ}_N(l, d)$ satisfies the 
recursion relation
\[
\widehat{AQ}_N(l, d) = \widehat{AQ}_{N - 1}(l, d) + \widehat{AQ}_{N - 1}(N - l, d -1).
\]
By the induction hypothesis, the above equation becomes
\begin{eqnarray*}
\widehat{AQ}_N(l, d) & = & { l - 1 \choose \lfloor \frac{d - 1}{2}\rfloor}
{N - l - 1 \choose \lfloor \frac{d}{2} \rfloor} +
{N - l - 1 \choose \lfloor \frac{d - 2 }{2} \rfloor} { l - 1 \choose \lfloor \frac{d - 1}{2}\rfloor} \\
& & \\
& = & { l - 1 \choose \lfloor \frac{d - 1}{2}\rfloor}
{N - l \choose \lfloor \frac{d}{2} \rfloor},
\end{eqnarray*}
and this proves (\ref{altp-formula}). The second equality above comes from the recurrence
relation of the Pascal triangle.

Finally, we observe that
\[  \fl
AQ_N(l, d) = \widehat{AQ}_N(l, d), \ \ {\rm for} \ d \ {\rm even}, \quad {\rm and} \quad
AQ_N(l, d) = \widehat{AQ}_N(l + 1, d), \ \ {\rm for} \ d \ {\rm odd}.
\]
These relations, together with formula (\ref{altp-formula}), prove the proposition.

\epf

\noindent
If we insert the result of the above proposition in proposition \ref{FN and AQ}, we get
the following corollary:

\begin{cor} The power series of ${\cal F}_N[u_c](n)$ around $u = 0$ is given by
\begin{eqnarray*}  \fl
{\cal F}_N[u_c](n) = I & + & \sum_{k = 1}^{\lfloor \frac{N}{2} \rfloor} (- 1)^k
(\frac{u}{N})^{2 k} \sum_{l = 0}^{N - 1}
{ l - 1 \choose k - 1}
{N - l \choose k} \cdot E_{\delta}(- 2l, n) \\
& & \\
& + & \sum_{k = 0}^{\lfloor \frac{N + 1}{2} \rfloor} (- 1)^k (\frac{u}{N})^{2 k + 1}
\sum_{l = 0}^{N - 1}
{ l \choose k} {N - l - 1\choose k}
\cdot E_{\delta}(- 2l, n) \cdot J.
\end{eqnarray*}
\end{cor}

We shall now prepare the necessary tools for the proof of theorem \ref{Beta1}.
First, we shall consider the appropriate limit of $AQ_N(l, d)$ when
$N$ goes to infinity.
The subset $\widehat{D}^{disc}_d(l)$ of the discrete ordered simplex
\[\fl
\Delta^{disc}_d(N) = \{ (l_1, l_2, \ldots, l_d) \in (\mathbb{N} \cup \{0 \})^d;
N - 1 \geq l_1 > l_2 > \ldots > l_d \geq 0 \}
\]
with the edge of size $N$ is given by one equation. Its size $AQ_N(l, d)$ is therefore of the order $N^{d - 1}$.

\begin{lemma} Let $\lambda$ be a real number in $[0, 1]$ and let $\{l _N\}_{N \in \mathbb{N}}$ be a sequence of positive integers
such that   $\lambda_N < N$ and $\lim_{N \to \infty} \frac{l_N}{N} = \lambda$.
Then we have
\begin{equation}\fl 
\lim_{N \to \infty} \frac{d!}{N^{d - 1}} AQ_N(l_N, d) = \left \{
\begin{array}{cl} \frac{1}{B(\frac{d}{2}, \frac{d}{2} + 1)}
\ \lambda^{\frac{d}{2} - 1} (1 - \lambda)^{\frac{d}{2}}   = 
p_{\beta}(\lambda ; \frac{d}{2}, \frac{d}{2} + 1)\ ; &  d \ \ {\rm even} \\
& \\
\frac{1}{B(\frac{d +1}{2}, \frac{d + 1}{2})} \
\lambda^{\frac{d - 1}{2}}  (1 - \lambda)^{\frac{d - 1}{2}}  = 
p_{\beta}(\lambda; \frac{d + 1}{2}, \frac{d + 1}{2})
\ ; & d \ \ {\rm odd}\ \ .
\end{array}
\right.
\label{AQlimit}
\end{equation}
\label{LBeta-1}
\end{lemma}

\noindent{\bf Proof:} We shall prove the formula only for even $d$. The proof for odd $d$
is essentially the same. For $ d = 2 m$, formula (\ref{AQexplicit}) gives
\[
AQ_N( l_N, d) = {l_N - 1\choose m - 1}{N - l_N \choose m}.
\]
This expression can be expanded into
\[
AQ_N(l_N, d) = \frac{1}{(m - 1)! \, m!} \prod _{k = 0}^{m - 2} ( (l_N - 1) - k)
\prod_{k = 0}^{m - 1} ((N - l_N) - k),
\]
which gives
\[ \fl
AQ_N( l_N, d) = \frac{1}{(m - 1)! \, m!}
\Bigl( (N \frac{l_N}{N} - 1)^{m - 1} + {\cal O}((N \lambda)^{m - 2} )\Bigr)
\Bigl( (N - N \frac{l_N}{N})^m + {\cal O}((N \lambda)^{m - 1}) \Bigr),
\]
and, due to $d - 1 = m + (m - 1)$ and $\lim_{N \to \infty}\frac{l_N}{N} = \lambda$,
\[
\lim _{N \to \infty} \frac{1}{N^{d - 1}} \,AQ_N(\lambda_N, d) = \frac{1}{(m - 1)! \, m!} \,
\lambda^{m - 1} \,(1 - \lambda)^m.
\]
The definition of the Euler beta function for positive integers gives
$
\frac{d!}{(m - 1)! \, m!} = \frac{1}{B(m, m + 1)},
$
and this proves formula (\ref{AQlimit}) for even $d$.

\epf

\noindent{\bf Proof of theorem \ref{Beta1}:}
{Recall the set $\widehat{D}_d^{disc}(l)$, given by formula (\ref{Ddisc}). 
Rescaling it by the factor $1/N$ gives the set
\[ \fl
\widehat{D}^{disc}_d(\frac{l}{N}) = \{(\frac{l_1}{N},  \frac{l_2}{N}\ldots , \frac{l_d}{N});  
\frac{N - 1}{N} \geq  \frac{l_1}{N} >  \frac{l_2}{N}  > \ldots  > \frac{l_d}{N} \geq 0, 
\sum_{j = 1}^d (- 1)^{j  - 1} l_j = l \}
\]
which contains the same number of points as $\widehat{D}_d^{disc}(l)$, but lies in the polytope  
$\widehat{D}_d(l)$. 
Let $D_d^{disc}(l)$ denote the orthogonal projection of 
$\widehat{D}_d^{disc}(l)$ on the hyperplane $\{(x_1, \ldots, x_{d - 1}, 0)\} \subset \rf^d$.
The number $\sharp \widehat{D}_d^{disc}(\frac{l}{N})$ of points in 
$\widehat{D}_d^{disc}(\frac{l}{N})$ is clearly equal to the number
of points $\sharp D_d^{disc}(\frac{l}{N})$ in the projection.
So, on the one hand, the number $\sharp D_d^{disc}(\frac{l}{N})$ is equal to $AQ_N(l; d)$,
while on the other, the value 
$\frac{1}{N^{d - 1}}\sharp D_d^{disc}(\frac{l}{N})$ is approximately equal
to the volume ${\rm Vol}(D_d(l))$ of the projection of $\widehat{D}_d(l)$ on the 
hyperplane $x_d = 0$ in $\rf^d$. 
Let now $\{\lambda _N\}_{N \in \mathbb{N}}$ 
be a sequence of rationals $l/N$ converging to $\lambda \in [0, 1]$.
We have
\[
\lim_{N \to \infty }\frac{1}{N^{d - 1}} \, AQ_N(N \lambda_N, d) = {\rm Vol}(D_d(\lambda)).
\]
The above expression, together with lemma \ref{LBeta-1}, proves our theorem
\ref{Beta1}.

\epf

\noindent {\bf Proof of theorem \ref{Beta2}:} The proof is an obvious adaptation of the proof of
lemma \ref{LBeta-1}. We only have to replace the particular values $m + 1$ and $m$ of the shape parameters
by an arbitrary pair $a$ and $b$ of integers. Then the same calculations as those performed
in the proof of lemma \ref{LBeta-1} yield the proof.

The number $c(N)^{- 1}$ is the discrete integral of the function
\[
Q_N(\lambda_N, a, b) = \frac{ (a + b + 1)!}{N^{a + b}} \
{ N \lambda_N - 1 \choose a } { N - N \lambda_N \choose b }
\]
over the discrete interval $\lambda_N \in (0, \frac{1}{N}, \ldots, \frac{N - 1}{N})$ with the volume form
$(1/N)$. Multiplying by $c(N)$ normalizes $Q_N$ to a probability mass function whose integral
has to be equal to $1$.
\epf

\section{General ordered alternating partitions}
\label{partitions}

In this section, our goal is to express the numbers $AP_N(l, d)$ of ordered alternating partitions
of $l$ into $d$ parts not greater than $N - 1$. The parts in $AP_N(l, d)$ {\it need not be
distinct}. The numbers $AP_N(l, d)$ will be expressed in terms of the transformation
${\cal G}_N$ rather than ${\cal F}_N$ which yielded the numbers $AQ_N(l, d)$.

Let $\vec{u} = (u_0, u_1, \ldots, u_{N - 1}) \in \rf ^N$
be a discrete function. Recall that ${\cal G}_N[\vec{u}]$ is given by
\[
{\cal G}_N[\vec{u}](n) = 
E_{\delta}(1, n) \cdot R(\frac{u_{N - 1}}{N}) \cdot E_{\delta}(1, n) 
\cdot R(\frac{u_{N - 2}}{N}) \cdots
E_{\delta}(1, n) \cdot R(\frac{u_0}{N}).
\]
We have seen in section \ref{Transforms} that
\[
{\cal G}_N[u](n) = \prod _{l = N - 1}^0 \Bigl( \cos {\frac{u_l}{N}}\cdot I + \sin{\frac{u_l}{N}}
\,
E_{\delta}(- 2 l , n) \cdot J \Bigr).
\]
Let us consider the derivative of order $d$, namely
$
D^d {\cal G}_N[\vec{u}](n)  =
(\frac{d}{d s})^d |_{s =0} \,{\cal G}_N[s \vec{u}](n).
$
The generalized Leibniz rule gives
\begin{equation}  \fl
D^d{\cal G}_N[\vec{u}](n) =
\sum_{\vec{k}} {d \choose k_0, k_1, \ldots k_{N - 1}} \! \! \!
\prod_{l = N - 1}^0 \Bigl(\
\cos{(s \frac{u_l}{N}) }\, I + \sin{(s \frac{u_l}{N}) } \,E_{\delta}(- 2 l, n) \cdot J \Bigr)^{(k_l)} \! \!
 |_{s = 0},
\label{F-1}
\end{equation}
where we sum over all $\vec{k} = (k_0, k_1, \ldots , k_{N - 1})$ such that
$
\| \vec{k} \|_1 = k_0 + k_1 + \ldots + k_{N - 1} = d.
$
Let
$
p(k) = k \Mod 2
$
be the parity of $k$, and
let us define the operator
\[
{\rm alt} \colon (\mathbb{N} \cup \{0\})^N \to \mathbb{N} \cup \{0\}
\]
by
\begin{eqnarray*}   \fl
{\rm alt}(\vec{k}) & = & \Bigl(\sum_{i = 1}^{k_{N - 1}} (- 1)^{i - 1} (N - 1) \Bigr) + \ldots +
(- 1)^{p(k_{N - 1}) + \ldots + p(k_{j + 1})}
\Bigl( \sum_{i = 1}^{k_j} (- 1)^{i - 1} j \Bigr) + ...\\
 \fl & + & (- 1)^{p(k_{N - 1}) + \ldots + p(k_{2})}
\Bigl( \sum_{i = 1}^{k_1} (- 1)^{ i - 1} 1 \Bigr).
\end{eqnarray*}
Another important function in this section is ${\rm odd}(\vec{k})$. By definition it is equal to
the number of odd components
in the integral vector $\vec{k} = (k_0, k_1, \ldots k_{N - 1}) \in (\mathbb{N} \cup \{0\})^N$.
The following two examples should clarify the formula for ${\rm alt}$.
Let first $N = 6$, $d = 4$, and
$
\vec{k} = (1, 1, 1, 1, 0, 0).
$
Then,
$
{\rm alt}(\vec{k}) = 3 - 2 + 1 - 0 = 2.
$
Here we get  ${\rm alt}(\vec{k}) = 2$ and ${\rm odd}(\vec{k}) = 4$.
Let now $N = 6$, $d = 8$, and
$
\vec{k} = (0, 1, 2, 2, 3, 0).
$
Then,
$
\rm{alt}(\vec{k}) = 4 - 4 + 4 - 3 + 3 - 2 + 2 - 1 = 3,
$
so   ${\rm alt}(\vec{k})= 3$ and ${\rm odd}(\vec{k}) = 2$. 
The values $\| \vec{k}\|_1 = d$
and ${\rm odd}(\vec{k})$ are equal only when all the components $k_j$ of $\vec{k}$ are
equal either to $0$ or to $1$.

The following proposition follows directly from the definition of the function ${\rm alt}$.
\begin{prop}
The function ${\rm alt}$ has the following three properties:
\begin{enumerate}
\item The even components $k_j$ of $\vec{k}$ do not contribute to ${\rm alt}(\vec{k})$.
\item Replacing any odd component $k_j$ of $\vec{k}$ by $1$ does not change ${\rm alt}$.
\item In the $L_1$-sphere $S^d_1 = \{\vec{k}; \| \vec{k}\|_1 = d \}$ we have
\[
\sum_{\vec{k}; \ {\rm alt}(\vec{k}) = l} \! \! \! \!1 = AP_N(l, d).
\]
\end{enumerate}
\label{alt}
\end{prop}

The factors of (\ref{F-1}) are equal to 
\[ \fl
(\frac{d}{ds})^{k_l} |_{s = 0} \Bigl( \cos{( s\, \frac{u_l}{N}) }\, I + 
\sin{(s \, \frac{u_l}{N})} E_{\delta}(- 2 l, n) 
\cdot J\Bigr)= (\frac{u_l}{N} )^{k_l} \,
E_{\delta}(-2 l , n) ^{p(k_l)} J^{k_l}.
\]
The exponential factor $E_{\delta}(- 2 l, n )$ appears if and only if $k_l$ is odd. 
Formula (\ref{F-1}) can be rewritten as
\[  \fl
D^d{\cal G}_N[\vec{u}](n) =
\sum_{m = 0}^d    (\frac{1}{N})^d \Bigl( \! \! \!
\mathop{\sum _{\|\vec{k}\|_1 = d}}_ {{\rm odd}(\vec{k})= m} \! \! \!
{d \choose k_0, k_1, \ldots, k_{N - 1}} u_0^{k_0} u_1^{k_1} \cdots u_{N - 1}^{k_{N - 1}}
\! \!\! \!
\mathop{\prod_{l_m > \ldots > l_1}}_{k_{l_j} {\rm odd} } \! \! \!
E_{\delta}(- 2 l_j , n) J^{k_{l_j}} \Bigr).
\]
The indices $l_j \in \{l_1, \ldots l_m\} \subset \{0, 1, \ldots, N- 1\}$
in the product at the end of our formula are those for which the component $k_{l_j}$
of the vector $\vec{k} = (k_0, k_1, \ldots , k_{N - 1}) $ is an odd integer. As in the
previous section, we use the relation
$
J \cdot E(l, n) = E(- l, n) \cdot J
$
and get
\begin{equation}  \fl
D^d{\cal G}_N[\vec{u}](n)=
\Bigl(
\sum _{\|\vec{k}\|_1 = d} (\frac{1}{N})^d
{d \choose k_0, k_1 \ldots, k_{N - 1}} u_0^{k_0} u_1^{k_1} \cdots u_{N - 1}^{k_{N - 1}} \cdot
E_{\delta}(- 2 \,{\rm alt}(\vec{k}), n )\Bigr) \cdot J^d.
\label{F-2}
\end{equation}
The first two properties of ${\rm alt}$ from proposition (\ref{alt}) allow us to replace the alternating sums
\[
l_m - l_{m - 1} + l_{m - 2} - \ldots + (- 1)^{m - 1} l_1, \quad m = 1, \ldots , d,
\]
where $l_j$ are odd,
by the values ${\rm alt}(\vec{k})$ of vectors $\vec{k}$ appearing in the sum in (\ref{F-2}).

\noindent {\bf Proof of proposition \ref{multinomial}:}
Let  $P_{N, d}$ denote the probability mass function of the binomial distribution,
\[
P_{N, d}(\vec{u}, \vec{k}) =
{d \choose k_0, k_1 \ldots, k_{N - 1}} u_0^{k_0} u_1^{k_1} \cdots u_{N - 1}^{k_{N - 1}}.
\]
We divide the set of non-negative integer valued functions
$\{\vec{k}\} \subset (\mathbb{N} \cup \{0\})^{N - 1}$ into disjoint subsets 
with respect to the values of the
function ${\rm alt}\colon \{\vec{k}\} \to \mathbb{N} \cup \{0\}$. This gives
\begin{equation}
D^d{\cal G}_N[\vec{u}](n) = \sum_{l = 0}^{N - 1} (\frac{1}{N})^d
\Bigl(
\mathop{\sum _{\vec{k}}}_ {{\rm alt}(\vec{k}) = l}
\! \! \! \! \!
P_{N, d}(\vec{u}, \vec{k}) \Bigr) \cdot E_{\delta}( - 2 l , n) \cdot J^d.
\label{F-3}
\end{equation}
We can now apply the discrete inverse linear Fourier transform and obtain
\begin{equation}
\sum _{\vec{k}; \ {\rm alt}(\vec{k}) = l}
P_{N, d}(\vec{u}, \vec{k}) \cdot I=   N^d
\Bigl( \sum_{n = 0}^{N - 1} 
E_{\delta}(2 l , n) \cdot D^d{\cal G}_N[\vec{u}](n) \Bigr) \cdot J^{ - d}.
\label{IF-1}
\end{equation}
This is the formula that we had to prove.

\epf

\noindent {\bf Proof of proposition \ref{APexpression}:} 
 Recall the differential operator ${\cal D}_d$, 
defined in proposition  \ref{APexpression} of the introduction. 
First, we observe that
\[
{\cal D}_d \Bigl(\sum _{\vec{k}; \ {\rm alt}(\vec{k}) 
= l} P_{N, d} (\vec{u}, \vec{k}) \Bigr) |_{\vec{u} = 0}
= d! \ AP_N(l, d).
\]
This follows from the third part of proposition \ref{alt} and from 
the fact that for two vectors $\vec{k}_a$ and $\vec{k}_b$
we have the following possibilities:
\[ \fl
\partial _{\vec{k}_a} \Bigl[{d \choose k^b_0, k^b_1 \ldots, k^b_{N - 1}}
u_0^{k^b_0} u_1^{k^b_1} \cdots u_{N - 1}^{k^b_{N - 1}}\Bigr]
=
\left \{\begin{array}{cl} p(u_0, \ldots , u_{N - 1}) \ ; & \|\vec{k}_a\|_1 < d \\
& \\
d! \, \delta _{a, b} \ ; & \| \vec{k} _a\|_1 = d \\
& \\
0 \; & \| \vec{k}_a\|_1 > d
\end{array}
\right. ,
\]
where $p(u_0, \ldots, u_{N - 1})$ is a polynomial without the constant term and
$\delta_{a, b}$ is the Kronecker delta. From (\ref{IF-1}) we now get
\[
AP_N(l, d) = \frac{N^d}{d!}\  {\cal D}_d
\Bigl( \sum_{n = 0}^{N - 1} 
E_{\delta}(2 l , n) \cdot D^d{\cal G}_N[\vec{u}](n) \Bigr) |_{\vec{u} = 0} \cdot J^{ - d},
\]
which proves the proposition.

\epf

\section{Acknowledgement}
I am grateful to Matja\v z Konvalinka for a very helpful discussion.
The research for this paper was supported in part by the research programme Analysis and
Geometry, P1- 0291, and the research project  Analysis, Equations, and Partial Differential Equations,
J1-9104,  funded by the Slovenian Research Agency.

\end{document}